\documentclass[12pt]{amsart}
\usepackage{amsmath}
\usepackage{amssymb}
\usepackage{amsmath,amssymb,amsthm,amscd}
\usepackage{geometry}
\usepackage{txfonts}
\usepackage{amsbsy}
\usepackage{graphicx,color}

\numberwithin{equation}{section}

\newtheorem{theo}{Theorem}[section]
\newtheorem{lemm}{Lemma}[section]
\newtheorem{coro}{Corollary}[section]

\def\begeq{\begin{equation}}
\def\endeq{\end{equation}}

\begin{document}

\title{Zero Comparison Theorem of Sturm-Liouville Type for Harmonic Heat Flow}
\author{Shi-Zhong Du}
\thanks{$^\dagger$Research partially supported by STU Scientific Research Foundation for Talents (SRFT-NTF16006), and partially supported by ``Natural Science Foundation of Guangdong Province" (2019A1515010605)}
\address{Department of Mathematics,
            Shantou University, Shantou, 51506, P. R. China.} \email{szdu@stu.edu.cn}

\renewcommand{\subjclassname}{%
  \textup{2010} Mathematics Subject Classification}
\date{Jun. 2020}

\begin{abstract}
   In this paper, we will prove a zero comparison theorem of Sturm-Liouville type for linearized harmonic heat flow in form of
       $$
            \frac{\partial v}{\partial t}=v_{rr}+\frac{m-1}{r}v_r-\frac{b(r,t)}{r^2}v, \ 0<r<R, t\in(t_1,t_2),
       $$
 where $b(r,t)$ is a bounded function satisfying
  $$
    b(0,t)\equiv m-1, \ b_r(0,t)\equiv 0,\ \ \forall t\in(t_1,t_2)
  $$
in case $v(0,t)\equiv0, \forall t\in(t_1,t_2)$.
\end{abstract}

\maketitle\markboth{Harmonic heat flow}{Zero comparison}

\vspace{20pt}

\section{Introduction}

A first zero comparison theorem was investigated by Nickel \cite{N} in 1962, and revived by Matano \cite{M} in 1982 and later by \cite{H}, \cite{A1}, \cite{A2} for one dimensional parabolic equation
   \begin{equation}\label{e1.1}
     u_t=a(x,t)u_{xx}+b(x,t)u_x+c(x,t)u, \ \ \forall x\in I, t\in(T_1,T_2),
   \end{equation}
where $I$ is an finite or infinite open interval in ${\mathbb{R}}$. They proved that the counting number ${\mathcal{Z}}(u(\cdot,t))$ of the zero set is monotone nonincreasing in time, together with some other properties. Properties of these type were later developed by Filippas-Kohn \cite{FK}, Chen-Pol\'{a}\v{c}ik \cite{CP} and Matano-merle \cite{MM} to
   \begin{equation}\label{e1.2}
     u_t=u_{rr}+\frac{n-1}{r}u_r+c(r,t)u, \ \ \forall r\in I, t\in(T_1,T_2)
   \end{equation}
in studying of semilinear heat equation
  \begin{equation}\label{e1.3}
    u_t=\triangle u+f(u), \ \ \forall x\in\Omega, t\in(T_1,T_2)
  \end{equation}
 of Fujita type.

  In this paper, we are concerning the equation of the form
  \begin{equation}\label{e1.4}
    \frac{\partial v}{\partial t}=v_{rr}+\frac{m-1}{r}v_r-\frac{b(r,t)}{r^2}v, \ 0<r<R, t\in(t_1,t_2),
  \end{equation}
where $b(r,t)$ is a bounded function satisfying
  $$
    b(0,t)\equiv m-1, \ b_r(0,t)\equiv 0,\ \ \forall t\in(t_1,t_2)
  $$
in case $v(0,t)\equiv0, \forall t\in(t_1,t_2)$. Moreover, the following boundary condition
  \begin{equation}\label{e1.5}
    v(r,t)\equiv0 \mbox{ or } v(r,t)\not=0,\ \ \ \forall r=0,R,\ \  t\in(t_1,t_2)
  \end{equation}
is imposed.  Equation \eqref{e1.4} comes from the linearization procedure of
   \begin{equation}\label{e1.6}
     \theta_t=\theta_{rr}+\frac{m-1}{r}\theta_r-\frac{m-1}{r^2}\sin \theta\cos \theta,\ \ \forall r\in I, t>0,
   \end{equation}
 which is satisfied by rotational symmetric solution $u\in C^2({\mathbb{R}}^m,{\mathbb{S}}^m), \ {\mathbb{S}}^m\subset{\mathbb{R}}^{m+1}$
   \begin{equation}\label{e1.7}
   u(x,t)=(u_1,u_2,\cdots, u_{m+1})=\Bigg(\frac{x}{|x|}\sin \theta(r,t),\cos \theta(r,t)\Bigg)
  \end{equation}
 of harmonic heat flow
     \begin{equation}\label{e1.8}
       \frac{\partial u}{\partial t}=\triangle_g u-A_u(\nabla u,\nabla u),
    \end{equation}
 where $g$ is the induced metric of ${\mathbb{S}}^m$ and $A_u: T_u{\mathbb{S}}^m\times T_u{\mathbb{S}}^m\to(T_u{\mathbb{S}}^m)^\perp$ is the second fundamental form of ${\mathbb{S}}^m\subset R^{m+1}$ at $u$. The readers may refer to the work of Coron-Ghidaglia \cite{CG} for the calculations. (see also \cite{CDY})

 The main purpose of this paper is to prove the following zero comparison theorem for \eqref{e1.4}.

\begin{theo}\label{t1.1}
  Let $v$ be a classical solution of \eqref{e1.4} on $[0,R]\times(t_1,t_2)$ which is not identical to zero and satisfies \eqref{e1.5} for some $0<R<+\infty$. We define
    $$
     {\mathcal{Z}}(v(\cdot,t))\equiv\sharp\Big\{r\in[0,R]\Big|\ v(r,t)=0\Big\}
    $$
  to be the zero number of $v(\cdot,t)$ counting multiplicity. Then

(i) ${\mathcal{Z}}(v(\cdot,t))<\infty$ for any $t_1<t<t_2$,

(ii) ${\mathcal{Z}}(v(\cdot,t))$ is a monotone non-increasing function in time $t$,

(iii) if $v(r_0,t_0)=v_r(r_0,t_0)=0$ for some $0\leq r_0\leq R$ and $t_1<t_0<t_2$, then
    $$
     {\mathcal{Z}}(v(\cdot,t))>{\mathcal{Z}}(v(\cdot,s))\ \ \mbox{ for any } t_1<t<t_0<s<t_2.
    $$
\end{theo}

\vspace{10pt}

We will firstly present some crucial lemmas in Section 2 and then complete the proof of the main theorem in Section 3.

\vspace{40pt}

\section{Crucial lemmas}

Comparing to equation in form of \eqref{e1.1} or \eqref{e1.2}, there is an extra singular term $\frac{b(r,t)}{r^2}v$ in \eqref{e1.4}. However, due to our assumption on
   $$
      b(0,t)\equiv m-1,\ \ b_r(0,t)\equiv0, \ \ \forall t\in(t_1,t_2),
   $$
the singular terms
   $$
    \frac{m-1}{r}v_r-\frac{b(r,t)}{r^2}v
   $$
are regular in family of functions
   $$
     v\in C^\infty([0,R)), \ \ v(0)=0.
   $$
So, we surely can prove a similar zero comparison theorem of Sturm-Liouville type. To show the main result, we need several crucial lemmas. The first one is the following.

\begin{lemm}\label{l2.1}
  Let $v$ be a solution to \eqref{e1.4} under boundary condition \eqref{e1.5} for
     $$
       v(0,t)\equiv0, \ \ \forall t\in(t_1,t_2).
     $$
  If $v$ is not identical to zero, then

  (i) for any $t_0\in(t_1,t_2)$, there exist nonnegative interger $k$ and a solution $\psi(r)$ of
    \begin{equation}\label{e2.1}
     \begin{cases}
      \psi_{rr}+\frac{m-1}{r}\psi_r-\frac{m-1}{r^2}\psi-\frac{1}{2}r\psi_r+\frac{k}{2}\psi=0, & \forall r>0\\
      \psi(0)=0, \psi_r(0)\not=0,
      \end{cases}
    \end{equation}
  such that $\varepsilon^{-k}v(\varepsilon r,t_0-\varepsilon^2)$ tends to $\psi(r)$ local uniformly in $C^{2+\alpha}_{loc}([0,+\infty)), \alpha\in(0,1)$ as $\varepsilon\downarrow0$, and

  (ii) for any $t_0\in(t_1,t_2)$ and $r_0\not=0$, there exist nonnegative interger $k$ and a solution $\psi(r)$ of
    \begin{equation}\label{e2.2}
     \begin{cases}
      \psi_{rr}-\frac{1}{2}r\psi_r+\frac{k}{2}\psi=0, & \forall r>0\\
      \psi(0)=0, \psi_r(0)\not=0,
      \end{cases}
    \end{equation}
  such that $\varepsilon^{-k}v(r_0+\varepsilon r,t_0-\varepsilon^2)$ tends to $\psi(r)$ local uniformly in $C^{2+\alpha}_{loc}([0,+\infty)), \alpha\in(0,1)$ as $\varepsilon\downarrow0$.
\end{lemm}

\noindent\textbf{Proof.} Away from origin, $\frac{b(r,t)}{r^2}$ is a bounded function. So, part (ii) is inferred from part (ii) of Lemma 2.5 in \cite{CP}. We remain to show that part (i) holds. In fat, by Lemma \ref{l2.2} below, there exists a nonnegative integer $k$, such that for some nonnegative integers $\alpha_0$ and $\beta_0$ satisfying
    $$
      \alpha_0+2\beta_0=k,
    $$
the following condition
   $$
      \frac{\partial^{\alpha_0+\beta_0}}{\partial^{\alpha_0} r\partial^{\beta_0} t}v(0,t_0)\not=0,\ \ \ \frac{\partial^{\alpha+\beta}}{\partial^\alpha r\partial^\beta t}v(0,t_0)=0
   $$
holds for all nonnegative integers $\alpha,\beta$ fulfilling $\alpha+2\beta<k$.

Setting
  $$
    u(\varepsilon,r)\equiv v(\varepsilon r,t_0-\varepsilon^2),
  $$
by Taylor's expansion, we have
   \begin{eqnarray*}
     u(\varepsilon,r)&=&\Sigma_{\alpha+2\beta=k}\frac{C^\alpha_k(-2)^\beta\beta!\frac{\partial^{\alpha+\beta}}{\partial r^\alpha\partial t^\beta}v(0,t_0)}{k!}r^\alpha\varepsilon^k\\
      &&+\Sigma_{\alpha+\beta+\gamma=k+1}\frac{C^\alpha_k(-2)^\beta}{k!}r^\alpha\frac{\partial^\gamma}{\partial\varepsilon^\gamma}(\varepsilon^\beta)\int^\varepsilon_0(\varepsilon-\tau)^k\frac{\partial^{\alpha+\beta}}{\partial r^\alpha\partial t^\beta}v(\tau r, t_0-\tau^2)d\tau.
   \end{eqnarray*}
When we denote
   $$
     \eta(r,\varepsilon)\equiv\varepsilon^{-k}v(\varepsilon r,t_0-\varepsilon^2),
   $$
there holds
   $$
     \Bigg(\Big|\frac{\partial^\beta\eta}{\partial\varepsilon^\beta}\Big|+ \Big|\frac{\partial^\beta\eta}{\partial r^\beta}\Big|\Bigg)(r,\varepsilon)\leq C_{\beta}R^{k+1+\beta},\ \ \forall r\in(0,R], \varepsilon\in(0,1]
   $$
for any $\beta\geq0$, as long as $R\geq1$. Moreover,
  $$
    \limsup_{\varepsilon\to0^+}|\eta(r,\varepsilon)|\not\equiv 0.
  $$

Therefore, if one writes $v$ into self-similar variables
   $$
     V(r,s)\equiv e^{\frac{k}{2}s}v(e^{-\frac{s}{2}}r, t_0-e^{-s}),
   $$
there hold
   \begin{equation}\label{e2.3}
     \Bigg(\Big|\frac{\partial^\beta V}{\partial s^\beta}\Big|+ \Big|\frac{\partial^\beta V}{\partial r^\beta}\Big|\Bigg)(r,\varepsilon)\leq C_{\beta}R^{k+\beta},\ \ \forall r\in(0,R], s\geq-\log t_0
   \end{equation}
and
   \begin{equation}\label{e2.4}
      \lim_{s\to+\infty}V(r,s)\not\equiv0.
   \end{equation}
Furthermore, $V$ satisfies
   \begin{equation}\label{e2.5}
     V_s-\Big(V_{rr}+\frac{m-1}{r}V_r\Big)+\frac{1}{2}rV_r-\frac{k}{2}V=-\frac{B(r,s)}{r^2}V,
   \end{equation}
where
   $$
    B(r,s)\equiv b(e^{-\frac{s}{2}}r, t_0-e^{-s}).
   $$
So, there exists a limiting function $V_\infty(r,s), r\in[0,\infty), s\geq0$ which is not identical to zero, such that
   $$
     V(r,s+l)\to V_\infty(r,s)\ \ \mbox{ uniformly in }  C^{2+\alpha}_{loc}([0,\infty)\times[0,+\infty))
   $$
holds as $l\to+\infty$. In order to show that $V_\infty(r,s)$ equals to a solution $\psi$ of \eqref{e2.1}, we need the following decaying estimation:

 For any $R\geq1$, let's take a cutoff function $\phi\in C^\infty([0,2R))$, such that
  $$
    \phi(r)=\begin{cases}
      1, & r\in[0,R/2],\\
      0, & r>R
    \end{cases}
  $$
and
   $$
    0\leq\phi(r)\leq1,\ |\phi'(r)|+|\phi''(r)|\leq 1, \ \ \forall r\geq0.
   $$
Setting
   $$
     \varphi(r,s)=\phi(e^{-\frac{s}{2}}r)
   $$
for $0\leq r\leq R_s\equiv Re^{\frac{s}{2}}R$, we have
   \begin{eqnarray}\label{e2.6}\nonumber
      & |\varphi_r|(r,s)\leq e^{-\frac{s}{2}}\chi_{[R_s/2,R_s]}, & |\varphi_{rr}|(r,s)\leq e^{-s}\chi_{[R_s/2,R_s]},\\
      & |\varphi_s|(r,s)\leq \chi_{[R_s/2,R_s]}, & \forall r\in[0,R_s], s\geq-\log t_0.
   \end{eqnarray}
Multiplying \eqref{e2.5} by $\rho r^{m-1}V_s\varphi^2, \rho(r)\equiv e^{-\frac{r^2}{4}}$ and integrating over $[0,R_s]$, one gets
   \begin{eqnarray*}
     \int^{R_s}_0\rho r^{m-1}V_s^2\varphi^2&=&-\int^{R_s}_0\rho V_r(V_sr^{m-1}\varphi^2)_r+(m-1)\int^{R_s}_0\rho r^{m-1}V_rV_s\varphi^2\\
     &&+\frac{k}{2}\int^{R_s}_0\rho r^{m-1}VV_s\varphi^2-\int^{R_s}_0\rho r^{m-3}B(r,s)VV_s\varphi^2
   \end{eqnarray*}
after integration by parts. Using \eqref{e2.6} and
   $$
     \Big|\frac{\partial B}{\partial s}\Big|=\Big|-\frac{1}{2}e^{-\frac{s}{2}}rb_r+e^{-s}b_t\Big|\leq Ce^{-\frac{s}{2}},
   $$
we have
  \begin{eqnarray}\label{e2.7}\nonumber
    &&\frac{d}{ds}\int^{R_s}_0\rho r^{m-1}\Bigg(\frac{1}{2}V_r^2-\frac{k}{4}V^2+\frac{B(r,s)}{2r^2}V^2\Bigg)\varphi^2=-\int^{R_s}_0\rho r^{m-1}V_s^2\varphi^2\\ \nonumber
   && \ \ \ \ \ \ \ \ \ \ \ \ \ \ \ \  \ \  \ \  \ \ \ \ \ \ \ \ \ \ \ \ \ \ \ \  \ \  \ \ +\frac{1}{2}\int^{R_s}_0\rho r^{m-3}\frac{\partial B}{\partial s}V^2\varphi^2-2\int^{R_s}_0\rho r^{m-1}V_rV_s\varphi\varphi_r\\ \nonumber
   && \ \ \ \ \ \ \ \ \ \ \ \ \ \ \ \  \ \  \ \  \ \ \ \ \ \ \ \ \ \ \ \ \ \ \ \  \ \  \ \ +2\int^{R_s}_0\rho r^{m-1}\Bigg(\frac{1}{2}V_r^2-\frac{k}{4}V^2+\frac{B(r,s)}{2r^2}V^2\Bigg)\varphi\varphi_s\\
   && \ \ \ \ \ \ \ \ \ \ \ \ \ \ \ \  \ \  \ \  \ \ \ \ \ \ \ \ \ \ \ \ \ \ \ \  \ \leq-\int^{R_s}_0\rho r^{m-1}V_s^2\varphi^2+Ce^{-\frac{s}{2}}.
  \end{eqnarray}
Noting that
   $$
     \int^{R_s}_0\rho r^{m-1}\Bigg(\frac{1}{2}V_r^2-\frac{k}{4}V^2+\frac{B(r,s)}{2r^2}V^2\Bigg)\varphi^2
   $$
is bounded from below by \eqref{e2.3}, it's infer from \eqref{e2.7} that
   \begin{equation}\label{e2.8}
     \int^\infty_{-\log t_0}\int^{R}_0\rho r^{m-1}V_s^2<+\infty
   \end{equation}
for any $R>0$. Consequently,
   $$
    \int_0^{\infty}\int_0^R\rho r^{m-1}|\partial_sV_\infty|^2\leq\liminf_{l\to+\infty}\int^\infty_l\int_0^R\rho r^{m-1}V_s^2=0,
   $$
and thus $V_\infty$ is a steady state of \eqref{e2.5}. Because of
  $$
    B(r,s)\to b(0,t_0)=m-1 \ \ \mbox{ as } s\to+\infty,
  $$
 to show that $V_\infty$ fulfils \eqref{e2.1}, we need only to verify that
    \begin{equation}\label{e2.9}
      V'_{\infty}(0)\not=0.
    \end{equation}
 Suppose on the contrary, then
    $$
      V_\infty(0)=V'_{\infty}(0)=0.
    $$
 Integrating \eqref{e2.1} over $[0,R]$, we have
   \begin{eqnarray*}
     V'_\infty(R)&=&-(m-1)\int^R_0\frac{V'_\infty(r)}{r}dt+(m-1)\int^R_0\frac{V_\infty(r)}{r^2}dr\\
       &&+\frac{1}{2}\int^R_0rV'_\infty(r)dr+\frac{k}{2}\int^R_0V_\infty(r)dr\\
       &=&-(m-1)\frac{V_\infty(R)}{R}+\frac{1}{2}RV_\infty(R)+\frac{k-1}{2}\int^R_0V_\infty(r)dr.
   \end{eqnarray*}
 So,
   $$
     \Big(r^{m-1}V_\infty(r)\Big)'-\frac{1}{2}r^mV_\infty(r)-\frac{k-1}{2}r^{m-1}\int^r_0V_\infty(r')dr'=0
   $$
for all $r>0$. Integrating again over $[0,R]$, we obtain
   \begin{eqnarray}\label{e2.10}\nonumber
     R^{m-1}V_\infty(R)&=& \frac{1}{2}\int^R_0r^mV_\infty(r)dr+\frac{k-1}{2}\int^R_0r^{m-1}\int^r_0V_\infty(r')dr'dr\\
      &=&\int^R_0\Bigg(\frac{1}{2}r^m+\frac{k-1}{2m}(R^m-r^m)\Bigg)V_\infty(r)dr.
   \end{eqnarray}
Next, we show that $V_\infty\equiv0$. If not, letting any $R_0\in\Big(0,\sqrt{\frac{m}{m+k-1}}\Big]$ such that
   $$
     |V_\infty(R_0)|=\max_{r\in[0,R_0]}|V_\infty(r)|,
   $$
it's deduced from \eqref{e2.10} that
   \begin{eqnarray*}
     R_0^{m-1}|V_\infty(R_0)|&\leq&|V_\infty(R_0)|\int^{R_0}_0\Big|\frac{1}{2}r^m+\frac{k-1}{2m}(R_0^m-r^m)\Big|dr\\
      &\leq&\frac{m+k-1}{2m}R_0^{m+1}|V_\infty(R_0)|\leq\frac{1}{2}R_0^{m-1}|V_\infty(R_0)|.
   \end{eqnarray*}
Therefore, one can conclude that $V_\infty(r)\equiv0$ in $\Big[0,\sqrt{\frac{m}{m+k-1}}\Big]$. And thus $V_\infty$ must be identical to zero on $[0,+\infty)$ by uniqueness of ordinary differential equation, which is contradicting with the non-triviality of $V_\infty$ mentioned above. So, \eqref{e2.9} holds and the lemma was now proven. $\Box$\\

We can now complete the proof of Lemma \ref{l2.1} by showing:

\begin{lemm}\label{l2.2}
   Under assumptions of Lemma \ref{l2.1}, the solution to \eqref{e1.4} can not vanish at $(0,t_0)$ with infinitely many order.
\end{lemm}

\noindent\textbf{Proof.} At first, after replacing $v$ by $r\widetilde{v}$, \eqref{e1.4} changes to
   \begin{equation}\label{e2.11}
     \widetilde{v}_t-\widetilde{v}_{rr}-\frac{m+1}{r}\widetilde{v}_r=-\frac{b(r,t)-(m-1)}{r^2}\widetilde{v}.
   \end{equation}
Suppose that $\widetilde{v}$ is not identical to zero, then there exist $\delta_0>0$ and $\delta_0<R_0<R-\delta_0$, such that
   $$
     |\widetilde{v}|(r,t)\geq\delta_0, \ \forall R_0-\delta_0\leq r\leq R_0+\delta_0, t_0-\delta_0<t<t_0+\delta_0.
   $$
Taking a cut-off function $\zeta\in C^\infty([0,+\infty))$ such that
   $$
     \zeta(r)=\begin{cases}
        1, &  0<r<R_0-\delta,\\
        0, & r>R_0+\delta
     \end{cases}
   $$
and $0<\zeta(r)<1, \forall r\in(R_0-\delta,R_0+\delta)$,  we have the function $\overline{v}(r,t)\equiv \zeta \widetilde{v}(r,t)+1-\zeta(r)$ satisfies that
  \begin{equation}\label{e2.12}
    \overline{v}_t=\overline{v}_{rr}+\frac{m+1}{r}\overline{v}_r-\overline{b}(r,t)\overline{v},
  \end{equation}
where
   \begin{eqnarray*}
     \overline{b}(r,t)&=&\frac{b(r,t)-(m-1)}{r^2}+\overline{v}^{-1}\Bigg\{\frac{b(r,t)-(m-1)}{r^2}(\zeta-1)\\
     &&+2\zeta_r\widetilde{v}_r+\zeta_{rr}\widetilde{v}+\frac{m-1}{r}\zeta_r\widetilde{v}-\zeta_{rr}-\frac{m-1}{r}\zeta_r\Bigg\}
   \end{eqnarray*}
is a bounded function on $[0,+\infty)\times(t_0-\delta_0,t_0+\delta_0)$

 Re-scaling
   $$
     w(r,s)=\overline{v}(e^{-\frac{s}{2}}r,t_0-e^{-s}),
   $$
we have $w$ satisfies that
  \begin{equation}\label{e2.13}
    w_s+Lw=-B(r,s)w,
  \end{equation}
where
   $$
     Lw\equiv-w_{rr}-\frac{m+1}{r}w_r+\frac{1}{2}rw_r=-(\rho r^{m+1})^{-1}\frac{d}{dr}\Big(\rho r^{m+1}w_r\Big)
   $$
and
   $$
     B(r,s)=e^{-s}\overline{b}(e^{-\frac{s}{2}}r,t_0-e^{-s}).
   $$
Considering the eigenvalue problem of unbounded self-adjoint operator $L$ on Hilbert space
    $$
      L^{2}_{\rho}([0,+\infty))\equiv\Bigg\{w\in L^2_{loc}((0,+\infty))\Big|\ \  \int^\infty_0\rho(r)r^{m+1}w^2(r)dr<\infty\Bigg\}
    $$
with $\rho(r)\equiv e^{-\frac{r^2}{4}}$ and inner product
   $$
    \Big<u,v\Big>_{L^{2}_{\rho}([0,+\infty))}\equiv\int^\infty_0\rho(r)r^{m+1}u(r)v(r)dr,
   $$
 it's not hard to show that the gaps of consecutive eigenvalues are bounded from below by $1/2$ (see for example \cite{S}). Noting also that $B(\cdot,s)\to0$ as $s\to+\infty$, the result in \cite{A1} (Lemma 5, Page 435) guarantees the decaying rate of $w$ can not be faster than exponential, and hence the conclusion of the lemma holds true. $\Box$\\

Considering the linear parabolic equation
  \begin{equation}\label{e2.14}
    v_t=a(x,t)v_{xx}+b(x,t)v_x+c(x,t)v, \ \ x\in(a,b), t\in(t_0,t_1)
  \end{equation}
with bounded coefficients, and imposing the boundary condition
   \begin{equation}\label{e2.15}
     v(a)\equiv0 \mbox{ or } \not=0, \ \ v(b)\equiv0 \mbox{ or } \not=0,\ \ \forall t\in(t_0,t_1),
   \end{equation}
Angenent has proven in \cite{A2} (see also Lemma 2.3 in \cite{CP}) the following lemma:

\begin{lemm}\label{l2.3}
  Letting $v$ be a classical solution to \eqref{e2.14} and \eqref{e2.15}, ${\mathcal{Z}}(v(\cdot,t))$ is finite and monotone non-increasing for any $t\in(t_0,t_1)$. Moreover, if $v_x(x_0,t^*)=0$ at $(x_0,t^*)\in(a,b)\times(t_0,t_1)$, we have
     \begin{equation}\label{e2.16}
      {\mathcal{Z}}(v(\cdot,t))<{\mathcal{Z}}(v(\cdot,s)), \ \ \forall t_0<s<t^*<t<t_1.
     \end{equation}
\end{lemm}

As a corollary of Lemam \ref{l2.1} and \ref{l2.3}, we obtain a parallel version to Lemma 2.6 in \cite{CP}:

\begin{lemm}\label{l2.4}
  Under the assumptions in Theorem \ref{t1.1}. For any $t_1<t_0<t_2$, there exists a small $\varepsilon_0>0$ such that

(i) ${\mathcal{Z}}(v(\cdot,t))<\infty$ for any $t\in(t_0-\varepsilon_0,t_0)$,

(ii) $v(\cdot, t)$ has only simple zeros for each $t\in(t_0-\varepsilon_0,t_0)$,

(iii) ${\mathcal{Z}}(v(\cdot,t))$ is constant as $t$ varies in $(t_0-\varepsilon_0,t_0)$.
\end{lemm}

\noindent\textbf{Proof.} By Lemma \ref{l2.1}, $v(\cdot,t)$ has only simple zero at origin for $t\in(t_0-\varepsilon_0,t_0)$ as long as $\varepsilon_0$ is chosen to be small. Being away from the origin, \eqref{e1.4} becomes a linear heat equation with bounded coefficients. So, one can utilize Lemma \ref{l2.3} to yield part (i). Part (ii) is a consequence of part (i) and Hopf's boundary lemma for non-degenerate parabolic equation. Finally, part (iii) follows from (i) and (ii). $\Box$\\

  Secondly, we can prove a similar version of Lemma 2.7 in \cite{CP} without difficulty:

\begin{lemm}\label{l2.5}
 Under the assumptions in Theorem \ref{t1.1}, for each $t_0\in(t_1,t_2)$ and $\varepsilon>0$ given in Lemma \ref{l2.4}, we write
   $$
    m+1\equiv\lim_{t\uparrow t_0}{\mathcal{Z}}(v(\cdot,t)).
   $$
Then the followings hold true:

(i) There are $m+1$ continuous functions
 $$
   \zeta_i: (t_0-\varepsilon_0,t_0)\to[0,1]\ (0\leq i\leq m)
 $$
such that
   \begin{eqnarray*}
     & 0\equiv\zeta_0(t)<\zeta_1(t)<\cdots<\zeta_{m-1}(t)<\zeta_m(t),&\\
     &\{r\in[0,1]|\ v(r,t)=0\}=\{\zeta_i(t)|\ 0\leq i\leq m\}.&
   \end{eqnarray*}

(ii) The limit $\zeta_i(t_0)\equiv\lim_{t\uparrow t_0}\zeta_i(t)$ exists for every $0\leq i\leq m$.

(iii) ${\mathcal{Z}}(v(\cdot,t_0))<\infty$; more precisely,
  $$
   \{r\in[0,1]|\ v(r,t_0)=0\}=\{\zeta_i(t_0)|\ 0\leq i\leq m\}.
  $$

(iv) $v(\cdot, t_0)$ has only simple zeros if and only if
  $$
   0=\zeta_0(t_0)<\zeta_1(t_0)<\cdots<\zeta_m(t_0).
  $$
\end{lemm}

\noindent\textbf{Proof.} Part (i) is clear a consequence of Lemma \ref{l2.4} and implicit function theorem. Suppose that (ii) is not true for some $0<i\leq m$, then
    $$
     r_-\equiv\liminf_{t\uparrow t_0}\zeta_i(t)<r_+\equiv\limsup_{t\uparrow t_0}\zeta_i(t).
    $$
Taking some $r_0\in(r_-,r_+)$, by part (ii) of Lemma \ref{l2.1}, there exist nonnegative integer $k$ and a solution $\psi(r)$ of \eqref{e2.3}, such that $\varepsilon^{-k}v(r_0+\varepsilon r,t_0-\varepsilon^2)$ tends to $\psi(r)$ local uniformly as $\varepsilon\downarrow0$. Now, letting $r_1>0$ such that $\psi(r_1)\not=0$, we have $v(r,t)\not=0$ for any $r=r_0+r_1\sqrt{t_0-t}, t\in(t_0-\varepsilon_0,t_0)$ and small positive number $\varepsilon_0$. However, by definition of $r_-$ and $r_+$, the nodal set of $v$ intersects with the curve $r=r_0+r_1\sqrt{t_0-t}, t\in(t_0-\varepsilon_0,t_0)$ infinitely many times. Contradiction holds and thus gives the proof of part (ii).

 To show the part (iii) we claim that if $\zeta_i(t_0)<\zeta_{i+1}(t_0)$ for some $0\leq i\leq m-1$, then $v(r,t_0)\not=0$ for $\zeta_i(t_0)<r\zeta_{i+1}(t_0)$. In fact, noting that $v(r,t)=0$ for $r=\zeta_i(t), \zeta_{i+1}(t)$ and $t\in(t_0-\varepsilon_0,t_0]$, the claim follows from the strong maximum principle of parabolic equation. So part (iii) holds true.

 Finally, the only if part of (iv) is a corollary of differential intermediate value theorem, while the if part follows from Hopf's boundary lemma for parabolic equation. The proof of Lemma \ref{l2.5} is completed. $\Box$\\

 Finally, a minor change of Lemma 2.8 in \cite{CP} yields the following lemma.

\begin{lemm}\label{l2.6}
 Under the assumptions of Theorem \ref{t1.1} and let $t_1<T_1<T_2<t_2$.

(i) If there are $m+1$ points $0=r_0<r_1<\cdots<r_m\leq R$ such that
   $$
    v(r_i,T_2)v(r_{i+1},T_2)<0\ \ (1\leq i\leq n),
   $$
then there also exist $m+1$ points $0=r'_0<r'_1<\cdots<r'_n\leq R$ such that
   $$
    v(r'_i,T_1)v(r'_{i+1},T_1)<0\ \ (1\leq i\leq n).
   $$

(ii) If $v(\cdot,T_2)$ has only simple zeros and ${\mathcal{Z}}(v(\cdot,T_2))=m+1$, then there are
   $$
    0=\zeta_0<\zeta_1<\cdots<\zeta_{m-1}<\zeta_m\leq R
   $$
such that $v(\zeta_i,T_1)=0$. In particular, ${\mathcal{Z}}(v(\cdot,T_1))\geq{\mathcal{Z}}(v(\cdot,T_2))$.
\end{lemm}

\noindent\textbf{Proof.} Letting $\Omega\subset[0,R]\times(t_1,t_2)$ be a connected component of support of $v$ containing $(\frac{r_i+r_{i+1}}{2},T_2)$, we claim first that for any $t\in(t_1,T_2)$, the line $[0,R]\times\{t\}$ intersects with $\Omega$. If not, then there exists $t_*\in(t_1,T_2)$, such that $\Omega$ lies above the line $(\frac{r_i+r_{i+1}}{2},t_*)$. Then, $\Omega$ forms a paraboloid type domain and $v$ vanishes on its boundary. By strong maximum principle, we conclude that $v\equiv0$ inside $\Omega$, which contradicts with the non-triviality of $v$. So, part (i) holds and part (ii) follows from intermediate value theorem. $\Box$\\

\vspace{40pt}

\section{Complete the proof of main theorem}

Now, we can complete the proof of Theorem \ref{t1.1} as follows.

\noindent\textbf{Proof of Theorem \ref{t1.1}:} In case $v(0,t)\not=0, \forall t\in(t_1,t_2)$, all zeros of $v$ located away from $r=0$, and hence conclusion follows from a slightly variant version of Lemma 2.2 in \cite{A1} (or Lemma 2.3 in \cite{CP}) since no singular term presence for $r\geq r_0>0$. In case $v(0,t)\equiv0, \forall t\in(t_1,t_2)$, if $v_r(0,t_0)$ is also vanishing, it yields from Lemma \ref{l2.1} that there exists $\varepsilon_0>0$, such that $r=0$ is only simple root of $v(\cdot,t)$ for $t\in(t_0-\varepsilon_0,t_0)$. So, part (i) of theorem \ref{t1.1} was shown in Lemma \ref{l2.5} (iii). To show parti (ii), fixing any $t_1<T_1<T_2<t_2$, one has ${\mathcal{Z}}(v(\cdot,s))\geq{\mathcal{Z}}(v(\cdot,T_2))$ and the zeros of $v(\cdot,s)$ are all simple for $s$ closing sufficiently to $T_2$ from below by Lemma \ref{l2.4} and 2.7'. Combining with Lemma \ref{l2.6}, we conclude that ${\mathcal{Z}}(v(\cdot,T_1))\geq{\mathcal{Z}}(v(\cdot,s))$ and obtain part (ii) of Theorem \ref{t1.1}. Finally, for arbitrary $t_1<T_1<T_0<T_2<t_2$, we want to show that $v(\cdot, T_0)$ has only simple zeros provided ${\mathcal{Z}}(v(\cdot,T_1))={\mathcal{Z}}(v(\cdot,T_2))=m+1$. Actually, by Lemma 2.6', if one takes $s$ to be closing sufficiently to $T_2$ from below, the zeros of $v(\cdot,s)$ are all simple. Moreover, ${\mathcal{Z}}(v(\cdot,s))=m+1$ by part (ii) of Theorem \ref{t1.1}. Thus, it's inferred from Lemma \ref{l2.6} (ii) that there exist
   $$
    0=\zeta_0<\zeta_1<\cdots<\zeta_{m-1}<\zeta_m=1
   $$
 such that $v(\zeta_i,T_0)=0$. So, all zeros of $v(\cdot,T_0)$ are simple and the proofs of Theorem \ref{t1.1} were completed. $\Box$\\

It's notable to remark that when the end point $r=R$ is replaced by a moving free boundary $r=R(t)$, conclusion of Theorem \ref{t1.1} still holds true. One need only using the transformation
   $$
     \overline{v}(r,t)=v(R^{-1}(t)r,t).
   $$
As a consequence of the theorem, we also have the following corollary:

\begin{coro}\label{c3.1}
  Let $v$ be a classical solution of \eqref{e1.4} on $[0,R]\times(t_1,t_2)$ or on $[0,+\infty)\times(t_1,t_2)$, which satisfies \eqref{e1.5} when $0<R<+\infty$. Suppose that for some $t_1<t_*<t^*<t_2$ and $r^*\in[0,R]$ (or $r^*\in[0,+\infty)$ respectively), there holds
     \begin{equation}\label{e3.1}
       v_r(r^*,t)=v(r^*,t)=0\ \ \forall t\in[t_*,t^*],
     \end{equation}
  then $v(r,t)\equiv0$.
\end{coro}

\noindent\textbf{Proof.} If $v$ satisfies \eqref{e1.5} for $0<R<+\infty$, conclusion follows from Lemma \ref{l2.1} since when $v$ not identical to zero, $z(v(\cdot,t))$ can drop only finitely many zeros and hence contradict with \eqref{e3.1}. In case $v$ is a solution of \eqref{e1.4} on $[0,+\infty)\times(t_1,t_2)$, a same reason can be applied to exclude the possibility of $|v|(r,t)>0$
when $(r,t)$ lies near some $(r_0,t_0)\in(r^*,+\infty)\times(t_*,t^*)$.
In fact, if not, then $v$ must be identical to zero in $[0,r_0)\times(t_*,t^*)$ by Theorem \ref{t1.1}. This contradicts with our assumption $|v|(r,t)>0$ near $(r_0,t_0)$. Therefore, $v(r,t)\equiv0$ for $r\geq0, t\in (t_*,t^*)$. Hence, it must be also identical to zero by Lemma \ref{l2.2}, since $v$ vanishes at some point $(0,t_0)$ with infinitely order. The proof was done. $\Box$\\

\vspace{40pt}

\section*{Acknowledgments}
The author would like to express his deepest gratitude to Professors Kai-Seng Chou and Xu-Jia Wang for their constant encouragements and warm-hearted helps.\\

\end{document}